\documentclass[12pt,leqno]{article}
\usepackage{amssymb,amsmath,amsfonts,amsbsy, xspace, latexsym, amscd}

\newtheorem{theorem}{Theorem}[section]

\newtheorem{corollary}{Corollary}[section]

\newtheorem{simpler notation}{Simpler Notation}[section]
\newtheorem{remark}{Remark}[section]
\newtheorem{example}{Example}[section]

\newcommand{\Proof}{\noindent \textbf{Proof:\;}}

\begin{document}


\title{Exponentiating $2\times 2$ and $3\times 3$ Matrices Done Right}

		\author{	Angel P. Popov (apopov@nws.aubg.bg)\\
			American University in Bulgaria\\
			2700 Blagoevgrad, Bulgaria\\
\and
			Todor D. Todorov (ttodotrov@calpoly.edu)\\
			Mathematics Department\\		
			California Polytechnic State University\\
		 San Luis Obispo, CA 93407, USA}
\date{}
\maketitle	

\begin{abstract} We derive explicit formulas for calculating $e^A$, $\cosh{A}$, $\sinh{A}, \cos{A}$
and
$\sin{A}$ for a given $2\times2$ matrix $A$. We also derive explicit formulas for $e^A$ for a given
$3\times3$ matrix $A$. These formulas are expressed exclusively in terms of the characteristic roots of $A$
and involve neither the eigenvectors of $A$, nor the transition matrix associated with a particular
canonical basis. We believe that our method has advantages  (especially if applied by non-mathematicians
or students) over the more conventional methods based on the choice of canonical bases. We support this point
with several examples for solving first order linear systems of ordinary differential equations with
constant coefficients. 
\end{abstract}

\section{Introduction} \label{S: Introduction} 

		The exponential $e^{A}$ of a square matrix $A$ and the related one-parameter family $e^{tA}$ are
important concepts in mathematics. Here is one example (among many others): Let
$\vec{x^\prime}(t)=A\vec{x}(t)$ be a system of first order homogeneous ordinary differential
equations with constant coefficients with initial conditions
$\vec{x}(0)=\vec{x}_0\in\mathbb{R}^n$, where $A$ is an $n\times n$ matrix with real entries. Then the
solution to the system is given by the formula $\vec{x}(t)=e^{tA}\,\vec{x}_0$ (Michael Artin~\cite{mArtin},
p. 140). 

	The exponential $e^{tA}$ is defined by the Taylor expansion:
$e^{tA}=\sum_{n=0}^\infty\frac{t^nA^n}{n!}$. For some particular matrices $A$ the exponential can be
easily calculated. For example, if $N$ is a nilpotent matrix of order $m$, i.e.
$N^m=O$, then 
$e^{tN}=\sum_{n=0}^{m-1}\frac{t^nN^n}{n!}$. If $P$ is a projection, i.e. $P^2=P$, then
$e^{tP}=I+\sum_{n=1}^\infty\frac{t^n P}{n!}=I-P+ e^tP$. If $A=\alpha I$ for some scalar $\alpha$,
then 
$e^{tA}=\sum_{n=0}^\infty\frac{t^n\alpha^n}{n!}I=e^{\alpha t}I$. From the Taylor
expansion it follows that if
$AB=BA$ for two square matrices, then $e^{A+B}=e^A\, e^B$. Also, if $A$ and $C$ are similar
matrices with a transition matrix $T$, i.e. $A=TC\,T^{-1}$, then
$e^{A}=Te^{C}\,T^{-1}$. 

	These facts can be found in many textbooks (Jerry Farlow et al~\cite{McDill}, p. 350)  and handbooks
(\cite{Handbook}, p. 132-133) on linear algebra. For a general matrix $A$, however, the rule for calculating
$e^{tA}$ becomes somewhat more complicated (A.I. Malcev~\cite{Malcev}, p. 118 or Michael
Artin~\cite{mArtin}, p. 480-482): We have to find the Jordan form
$T^{-1}AT={\rm diag}(C_1, C_2,\dots,C_k)$ of the matrix $A$, where $C_i$'s are Jordan cells. Thus the
formula
$e^{A}=Te^{C}\,T^{-1}$ becomes $e^{A}=T{\rm diag}(e^{C_1}, e^{C_2},\dots,e^{C_k})\,T^{-1}$. On the other
hand, each cell can be presented in the form
$C_i=\lambda_iI+N$, where $N^{m_i}=O$, $\lambda_i$ is the corresponding eigenvalue of $A$ and $m_i={\rm
size}(C_i)$. As a result, $e^{C_i}=e^{\lambda_i}\sum_{n=0}^{m_i-1}\frac{N^n}{n!}$. The procedure is
universal and it looks attractive but it has the following disadvantages: (a) The algorithm for calculating
the transition matrix $T$ and its inverse $T^{-1}$ is difficult  and time consuming. It is especially
complicated in the case of multiple characteristic roots, when the calculations of the generalized
eigenvectors require skills that are too advanced for students in a typical course on ordinary differential
equations and linear algebra. (b) The original framework
$\mathbb{R}^n$ must be extended to $\mathbb{C}^n$, which again might be confusing for students.  Our
observation is that students and non-mathematicians (the latter might have rusty knowledge of linear
algebra) who rarely calculate exponentials of matrices often have considerable difficulty calculating
$e^{tA}$ even for
$2\times 2$ and $3\times 3$ matrices. 

	In this article we derive explicit formulas for calculating $e^{tA}$ for a $2\times
2$ and $3\times 3$ matrix $A$ with real entries. In the case of $2\times2$ matrices $A$ we also derive
explicit formulas for $\cosh{A}$, $\sinh{A}, \cos{A}$ and
$\sin{A}$. These formulas are expressed exclusively in terms of the characteristic roots of $A$
and involve neither the eigenvectors of $A$, nor the transition matrix associated with a particular
canonical basis. We believe that our formulas are
suitable for handbooks in the sense that one can apply them with limited background in linear algebra, in
particular, without having slightest idea what a ``transition matrix is''. The formulas are derived with
the help the characteristic polynomial $f_A(\lambda)=\det(A-\lambda I)$ of $A$ and the Cayley-Hamilton
theorem which  says that $f_A(A)=0$. But, again, the formulas can be applied without knowledge of
the Cayley-Hamilton theorem. In Corollary~\ref{C: Eigenvectors, Canonical Forms and
Bases} and  Corollary~\ref{C: Eigenvectors, Canonical Forms and
Bases 3by3} we identify some matrices closely related to $A$ which allow easy
calculation of the invariant subspaces of $A$ and the canonical bases (if needed). We recommend our method
for teaching exponentials in a course on linear algebra and ordinary differential equations. To support this
point we present several examples for solving linear systems of ordinary differential equations.

	For other methods for calculating the exponential of a matrix not mentioned in our paper we refer to
the survey article (Cleve Moler and Charles Van Loan~\cite{MolerVanLoan}), where the reader will find more
references to the subject. 

	In what follows $I$ and $O$ denote the identity and zero square matrices,
respectively.  We denote by
$\left[\,\vec{x}\,|\,\vec{y}\,\,|\,\vec{z}\,|\, \dots\right]$ the matrix with column-vectors
$\vec{x},\, \vec{y},\, \vec{z},\dots\in\mathbb{R}^n$. Also,
$E_{\lambda=\mu}(A)$ denotes the eigenspace of $A$ corresponding to $\mu\in\mathbb{C}$. 

\section{Exponential of a $2\times2$ Matrix}\label{S: Exponential of a 2times2-Matrix}

	We derive formulas for calculating the exponential of a given $2\times 2$-matrix $A$ based exclusively on
the characteristic roots of $A$. Our formulas are easy to memorize and simple to use in the sense that: (a)
our framework is always $\mathbb{R}^2$ (never $\mathbb{C}^2$); (b) our formulas involve linear
operations between matrices only; (c) our formulas can be handled with limited background in linear algebra
- they involve neither the eigenvectors of $A$, nor a particular canonical basis (although we need the
characteristic roots of $A$).  We illustrate this point with several examples of linear systems of ordinary
differential equations. Although we do not need a change of the basis to calculate the exponential
$e^{tA}$, we also derive explicit formulas for a canonical basis and a transition matrix
$T$  for which  $C=T^{-1}AT$ is exceptionally simple for the purpose of the exponentiation
(Corollary~\ref{C: Eigenvectors, Canonical Forms and Bases}).

\begin{theorem}\label{T: Main 1} Let $A$ be a $2\times2$-matrix with real entries and let
$\lambda_1$ and $\lambda_2$ be the characteristic roots of $A$. 

		\textbf{Case 1}: If $\lambda_1=\lambda_2$ (real), then $A=\lambda_0 I+N$,
where $\lambda_0=\lambda_1=\lambda_2$ and the matrix $N=A-\lambda_0 I$ satisfies the equation $N^2=O$.
Consequently, for every real (or complex)
$t$ we have
\begin{equation}\label{E: Parabolic}
e^{tA}=e^{\lambda_0 t}(I+tN).
\end{equation}

		\textbf{Case 2}: If $\lambda_{1,2}=\alpha\pm i\omega$  for some $\alpha,\, 
\omega\in\mathbb{R},\, \omega \not=0$, then $A=\alpha I+\omega J$, where the matrix
$J=\frac{1}{\omega}(A-\alpha I)$ satisfies the equation $J^2=-I$. Consequently, for every real (or complex)
$t$ we have
\begin{equation}\label{E: Elliptic}
e^{tA}=e^{\alpha t}[(\cos{\omega t}) I+(\sin{\omega t})J].
\end{equation}

	\textbf{Case 3}: If $\lambda_1\not=\lambda_2$ (both real), then $A=\alpha I+\beta J$, where 
$\alpha=\frac{\lambda_1+\lambda_2}{2},\; \beta=\frac{\lambda_1-\lambda_2}{2}$
and the matrix $J=\frac{1}{\beta}(A-\alpha I)$ satisfies the equation $J^2=I$. Consequently, for every
real (or complex)
$t$ we have
\begin{equation}\label{E: Hyperbolic}
e^{tA}=e^{\alpha t}\{\,(\cosh{\beta t}) I+(\sinh{\beta
t})J\;\}\overset{or}=\frac{1}{2}\left[e^{\lambda_1t}(I+J)+e^{\lambda_2t}(I-J)\right].
\end{equation}
\end{theorem}
\Proof Case 1: We first show that $N^2=0$.
Indeed,
$N^2=(A-\lambda_0I)^2=A^2-2\lambda_0A+\lambda_0^2I=O$, by the Cayley-Hamilton
theorem, since $2\lambda_0={\rm tr}(A)$ and $\lambda_0^2=\det(A)$.  Thus we have
\[
e^{tA}=e^{t(\lambda_0 I+N)}=e^{\lambda_0t I}\cdot
e^{tN}=e^{\lambda_0t}I(I+tN+0+\dots)=e^{\lambda_0t}\,(I+tN),
\]
as required. 

	 Case 2: We first show that $J^2=-I$. Indeed, we have 
$A^2-2\alpha I+(\alpha^2+\omega^2)I=O$, by the Cayley-Hamilton
theorem, since $2\alpha=\lambda_1+\lambda_2={\rm tr}(A)$ and
$\alpha^2+\omega^2=\lambda_1\lambda_2=\det(A)$. Thus
$J^2=(\frac{1}{\omega}(A-\alpha I))^2=\frac{1}{\omega^2}(A^2-2\alpha
I+\alpha^2I)=\frac{1}{\omega^2}(A^2-2\gamma
I+(\alpha^2+\omega^2)I-\omega^2I)=\frac{1}{\omega^2}(O-\omega^2I)=-I$. Next, we calculate
\begin{align}\notag
&e^{tA}=e^{t(\alpha I+\omega J)}= e^{\alpha t I}\, e^{\omega tJ}=
e^{\alpha t}\left\{\sum_{n=0}^\infty\frac{(\omega t)^{2n}J^{2n}}{(2n)!}
+\sum_{n=0}^\infty\frac{(\omega t)^{2n+1}J^{2n+1}}{(2n+1)!}\right\}=\\\notag
&e^{\alpha t}\left\{\sum_{n=0}^\infty\frac{(-1)^n(\omega t)^{2n}}{(2n)!}I
+\sum_{n=0}^\infty\frac{(-1)^n(\omega t)^{2n+1}}{(2n+1)!}J\right\}=
e^{\alpha t}\left\{(\cos{\omega t}) I
+(\sin{\omega t})J\right\},
\end{align}
as required. 

	Case 3: We shall, first, show 
that $J^2=I$. Indeed, we have 
$A^2-2\alpha I+(\alpha^2-\beta^2)I=O$, by the Cayley-Hamilton
theorem, since
$2\alpha=\lambda_1+\lambda_2={\rm tr}(A)$ and
$\alpha^2-\beta^2=\lambda_1\lambda_2=\det(A)$. Thus
$J^2=(\frac{1}{\beta}(A-\alpha I))^2=\frac{1}{\beta^2}(A^2-2\alpha
I+\alpha^2I)=\frac{1}{\beta^2}(A^2-2\alpha
I+(\alpha^2-\beta^2)I+\beta^2I)=\frac{1}{\beta^2}(0+\beta^2I)=I$. Next, we calculate
\begin{align}\notag
&e^{tA}=e^{t(\alpha I+\beta J)}=e^{\alpha t I}\cdot e^{\beta t J}=e^{\alpha
t}\left\{\sum_{n=0}^\infty\frac{(\beta t)^{2n}J^{2n}}{(2n)!} +\sum_{n=0}^\infty\frac{(\beta
t)^{2n+1}J^{2n+1}}{(2n+1)!}\right\}=\\\notag &=e^{\alpha t}\left\{\sum_{n=0}^\infty\frac{(\beta
t)^{2n}}{(2n)!}I +\sum_{n=0}^\infty\frac{(\beta t)^{2n+1}}{(2n+1)!}J\right\}= e^{\alpha
t}\left\{(\cosh{\beta t}) I +(\sinh{\beta t})J\right\},
\end{align}
as required. 
$\blacktriangle$

	The formulas (\ref{E: Parabolic})-(\ref{E:
Hyperbolic}) show that we do not need the eigenvectors of $A$ in order to calculate $e^{tA}$. However, the
eigenvectors of $A$ as well as canonical forms and canonical bases (if needed for classification, graphing,
etc.) can be easily extracted from the matrices
$N$ and $J$. The next corollary follows easily from the above theorem and we leave the proof to the
reader. 
\begin{corollary}[Eigenvectors, Canonical Forms and Bases]\label{C: Eigenvectors, Canonical Forms and
Bases} Under the notation of the previous theorem we have the following:

	\textbf{Case 1:} If $N\not=O$, then $E_{\lambda=\lambda_0}(A)= {\rm Im}(N)$. In particular,
either of the non-zero columns of $N$ is an eigenvector of $A$. We
have
$A=T\begin{pmatrix}
\lambda_0& 1\\
0&\lambda_0
\end{pmatrix}T^{-1},$ where $T=[\,N\vec{x}\,|\,\vec{x}\,]$ and $\vec{x}\in
\mathbb{R}^2$ such that $N\vec{x}\not=\vec{0}$. If $N=O$, then $A=\lambda_0 I$.

	\textbf{Case 2:} The matrix
$A$ does not have eigenvectors in $\mathbb{R}^2$. We
have $A=T\begin{pmatrix}
\alpha & \omega\\
-\omega &\alpha
\end{pmatrix}
T^{-1}$, where $T=[\,J\vec{x}\,|\,\vec{x}\, ]$ and $\vec{x}\in\mathbb{R}^2, \vec{x}\not=\vec{0}$. Notice
that the matrix $\begin{pmatrix}
\alpha & \omega\\
-\omega &\alpha
\end{pmatrix}$ is conformal (i.e. the corresponding transformation preserves the angles in $\mathbb{R}^2$).
Alternatively, in 
$\mathbb{C}^2$ the matrix $A$ is diagonalizable with eigenvalues $\alpha\pm i\omega$ and
eigenvectors $(I\mp iJ)\vec{x},\, \vec{x}\not=\vec{0}$, respectively.

	\textbf{Case 3:} We have
$E_{\lambda=\lambda_{1,2}}(A)= {\rm Im}(I\pm J)$, respectively.  We
have $A= T\begin{pmatrix}
\lambda_1& 0\\
0&\lambda_2
\end{pmatrix}
T^{-1}$,
where $T=[\, (I+J)\vec{x}\,|\, (I-J)\vec{y}\,]$ with $\vec{x}, \vec{y}\in \mathbb{R}^2$ such that 
$(I+J)\vec{x}\not=\vec{0}$ and $(I- J)\vec{y}\not=\vec{0}$. 
\end{corollary}

	Here are several examples tested in class by the second author in a course on
linear algebra and ordinary differential equations. We strongly recommend these formulas for teaching
exponentials.

\begin{example}  We shall find the solution of the initial value problem without 
involving eigenvectors:
\begin{align}\notag
&x^\prime=3x+2y,\\\notag
&y^\prime=-8x-5y,\notag
\end{align}
with initial conditions $x(0)=1,\; y(0)=-1$. We have $A=\begin{pmatrix}3 & 2\\
																									-8 & -5\end{pmatrix}$ and $\vec{x}_0=\begin{pmatrix}1\\
																								                    	-1 \end{pmatrix}$.
The characteristic polynomial is $\lambda^2+2\lambda+1$ with roots 
$\lambda_1=\lambda_2=-1$ (Case 1). We calculate
\[
N=A+I=\begin{pmatrix}3 & 2\\
																									-8 & -5\end{pmatrix}+\begin{pmatrix}1 & 0\\
																									0 & 1\end{pmatrix}=\begin{pmatrix}4 & 2\\
																									-8 & -4\end{pmatrix}.
\]
We leave to the reader to verify that $N^2=O$. For the exponential of
$A$ we apply formula (\ref{E: Parabolic}):
\begin{align}\notag
&\vec{x}(t)=e^{tA}\, \begin{pmatrix}1\\
																								     -1 \end{pmatrix}= e^{-t}\begin{pmatrix}1+4t, &2t\\
																									-8t, & 1-4t\end{pmatrix}\begin{pmatrix}1\\
																								     -1 \end{pmatrix}=e^{-t}\begin{pmatrix}1+2t\\
																								    -1-4t
\end{pmatrix}.
\end{align}
Thus $x=e^{-t}(1+2t),\; y=-e^{-t}(1+4t)$.
\end{example}
\begin{remark} If a canonical form for $A$ and a transition matrix are still needed, we can
easily calculate them with the help of Corollary~\ref{C: Eigenvectors, Canonical Forms and
Bases}. For example, $T^{-1}AT=\begin{pmatrix}-1 & 1\\
																									0 & -1\end{pmatrix}$, where $T=[N\vec{x}\,|\, \vec{x}\,]=\begin{pmatrix}2 & 0\\
																									-4 & 1\end{pmatrix}$ for $\vec{x}=(0, 1)$. Notice that either of the columns
of $N$ is an eigenvector of $A$.
\end{remark}
\begin{example}  We shall find the solution to the initial value problem without 
involving eigenvectors:
\begin{align}\notag
&x^\prime=y,\\\notag
&y^\prime=-5x-2y,\notag
\end{align}
with initial conditions $x(0)=2,\; y(0)=1$. We have $A=\begin{pmatrix} 0& 1\\
																									-5 & -2\end{pmatrix}$ and $\vec{x}_0=\begin{pmatrix}2\\
																								                    	1 \end{pmatrix}$.
The characteristic polynomial is $\lambda^2+2\lambda+5$ with roots
$\lambda_{1,2}=-1\pm 2i$, thus $\alpha=-1$ and $\omega=2$ (Case 2). We calculate
\[
J=\frac{1}{\omega}(A-\alpha I)=\frac{1}{2}\left[\begin{pmatrix}0 & 1\\
																									-5 & -2\end{pmatrix}+\begin{pmatrix}1 & 0\\
																									0 & 1\end{pmatrix}\right]=\frac{1}{2}\begin{pmatrix}1 & 1\\
																									-5 & -1\end{pmatrix}.
\]
We leave to the reader to verify that $J^2=-I$. 
For the solution to the system we apply formula (\ref{E: Elliptic}):
\begin{align}\notag
&\vec{x}=e^{tA}\vec{x}_0=e^{-t}\left[\,(\cos{2t})
\begin{pmatrix}1 & 0\\
																									0 & 1\end{pmatrix}+(\sin{2t})\frac{1}{2}\begin{pmatrix}1 & 1\\
																									-5 & -1\end{pmatrix}\;\right]\begin{pmatrix}2\\
																								                    	1 \end{pmatrix}=\\\notag
&=e^{-t}\begin{pmatrix}\cos{2t}+\frac{1}{2}\sin{2t}, & \frac{1}{2}\sin{2t}\\
																									-\frac{5}{2}\sin{2t}, & \cos{2t}-\frac{1}{2}\sin{2t}\end{pmatrix}
\begin{pmatrix}2\\
																								                    	1
\end{pmatrix}=
e^{-t}\begin{pmatrix}2\cos{2t}+\frac{3}{2}\sin{2t}\\
																								                    	\cos{2t}-\frac{11}{2}\sin{2t}\end{pmatrix}
\end{align}
Thus $x=e^{-t}(2\cos{2t}+\frac{3}{2}\sin{2t}),\; y=e^{-t}(\cos{2t}-\frac{11}{2}\sin{2t})$.
\end{example}
\begin{remark} If a canonical form for $A$ and the corresponding transition matrix are still needed, we can
easily calculate them with the help of Corollary~\ref{C: Eigenvectors, Canonical Forms and
Bases}. For example, $T^{-1}AT=\begin{pmatrix}-1 & 2\\
																									-2 & -1\end{pmatrix}$, where $T=[J\vec{x}\,|\, \vec{x}\,]=\begin{pmatrix}1 & 2\\
																									-5 & 0\end{pmatrix}$ for $\vec{x}=(2, 0)$. In the framework of $\mathbb{C}^2$ we
have $T^{-1}AT=\begin{pmatrix}-1+2i & 0\\
																									0 & -1-2i\end{pmatrix}$, where
$T=[(I-iJ)\vec{x}\,|\,(I+iJ)\vec{x}\,]=\begin{pmatrix}2+i & 2-i\\
																									-5i & 5i\end{pmatrix}$ for $\vec{x}=(2, 0)$.
\end{remark}
\begin{example}  We shall find the solution to the initial value problem without 
involving eigenvectors:
\begin{align}\notag
&x^\prime=5x-y,\\\notag
&y^\prime=3x+y,\notag
\end{align}
with initial conditions $x(0)=1,\; y(0)=2$. We have $A=\begin{pmatrix}5 & -1\\
																									3 & 1\end{pmatrix}$ and $\vec{x}_0=\begin{pmatrix}1\\
																								                    	2 \end{pmatrix}$.
The characteristic polynomial is $\lambda^2-6\lambda+8$ with roots
$\lambda_1=4, \lambda_2=2$ (Case 3). We calculate $\alpha= 3,\, \beta=1$ and 
\[
J=\frac{1}{\beta}(A-\alpha I)=\begin{pmatrix}5 & -1\\
																									3 & 1\end{pmatrix}-\begin{pmatrix}3 & 0\\
																									0 & 3\end{pmatrix}=\begin{pmatrix}2 & -1\\
																									3 & -2\end{pmatrix}.
\]
We leave to the reader to verify that $J^2=I$. 
For the solution to the system we apply formula (\ref{E: Hyperbolic}):
\begin{align}\notag
&\vec{x}(t)=e^{tA}\, \begin{pmatrix}1\\
																								     2 \end{pmatrix}= e^{3t}\begin{pmatrix}\cosh{t}+2\sinh{t}, &
-\sinh{t}\\
																									3\sinh{t}, & \cosh{t}-2\sinh{t}\end{pmatrix}\begin{pmatrix}1\\
																								     2 \end{pmatrix}=\\\notag
&=e^{3t}\begin{pmatrix}\cosh{t}\\
																								    2\cosh{t}-\sinh{t}
\end{pmatrix}.
\end{align}
Thus $x=e^{3t}\cosh{t}\overset{or}=\frac{1}{2}(e^{4t}+e^{2t})$ and 
$y=e^{3t}(2\cosh{t}-\sinh{t})\overset{or}=\frac{1}{2}(e^{4t}+3e^{2t})$.
\end{example}
\begin{remark} If a canonical form for $A$ and the transition matrix are still needed, we can
easily calculate them with the help of Corollary~\ref{C: Eigenvectors, Canonical Forms and
Bases}. For example, $T^{-1}AT=\begin{pmatrix}4 & 0\\
																									0 & 2\end{pmatrix}$, where $T=[(I+J)\vec{x}\,|\, (I-J)\vec{x}\,]=\begin{pmatrix}-1
& 1\\-1 & 3\end{pmatrix}$ for $\vec{x}=(0, 1)$. Notice that the columns of $T$ are eigenvectors of $A$.
\end{remark}
\section{Some Applications}

	If $A$ is a square matrix, we define the matrices  $\cosh{A}, \sinh{iA},
e^{iA}, \cos{A}$ and $\sin{A}$ by the corresponding Taylor series. For
example, $\cosh{A}=\sum_{n=0}^\infty\frac{A^{2n}}{(2n)!}$, and the rest are defined similarly. In this
section we derive explicit formulas for these elementary functions which follow easily from Theorem~\ref{T:
Main 1}.
\begin{corollary}[Hyperbolic Functions]\label{C: Sinh} Let $A$ be a $2\times2$-matrix with real entries
and let
$\lambda_1$ and $\lambda_2$ be the characteristic roots of $A$. 

		\textbf{Case 1}: If $\lambda_1=\lambda_2$ (real), then 
\begin{equation}\label{E: ParabolicSin}\notag
\cosh{A}=(\cosh{\lambda_0})I+(\sinh{\lambda_0})N,\quad \sinh{A}=(\sinh{\lambda_0})I+(\cosh{\lambda_0})N, 
\end{equation}
where $\lambda_0=\lambda_1=\lambda_2$ and $N=A-\lambda_0 I$.

		\textbf{Case 2}: If $\lambda_{1,2}=\alpha\pm i\omega$  for some $\alpha,\, 
\omega\in\mathbb{R},\, \omega \not=0$, then
\begin{align}\label{E: EllipticSin}\notag
&\cosh{A}=(\cosh{\alpha})(\cos{\omega})I+(\sinh{\alpha})(\sin{\omega})J,\\\notag
&\sinh{A}=(\sinh{\alpha})(\cos{\omega})I+(\cosh{\alpha})(\sin{\omega})J,
\end{align}
where $J=\frac{1}{\omega}(A-\alpha I)$. 

	\textbf{Case 3}: If $\lambda_1\not=\lambda_2$ (both real), then 
\begin{align}\notag
&\cosh{A}=(\cosh{\alpha})(\cosh{\beta})I+(\sinh{\alpha})(\sinh{\beta})J,\\\notag
&\sinh{A}=(\sinh{\alpha})(\cosh{\beta})I+(\cosh{\alpha})(\sinh{\beta})J,
\end{align}
where $\alpha=\frac{\lambda_1+\lambda_2}{2},\; \beta=\frac{\lambda_1-\lambda_2}{2}$
and $J=\frac{1}{\beta}(A-\alpha I)$.
\end{corollary}
\Proof The results follow directly from Theorem~\ref{T: Main 1} and the formulas
$\cosh{x}=\frac{1}{2}(e^x+e^{-x})$ and $\sinh{x}=\frac{1}{2}(e^x-e^{-x})$. $\blacktriangle$
\begin{corollary}[Complex Exponent]\label{C: Complex Exponent} Let $A$ be a $2\times2$-matrix with real
entries and let
$\lambda_1$ and $\lambda_2$ be the characteristic roots of $A$. 

		\textbf{Case 1}: If $\lambda_1=\lambda_2$ (real), then 
\begin{equation}\label{E: ParabolicC}\notag
e^{\pm iA}=e^{\pm i\lambda_0}(I\pm iN),
\end{equation}
where
$\lambda_0=\lambda_1=\lambda_2$ and $N=A-\lambda_0 I$.

		\textbf{Case 2}: If $\lambda_{1,2}=\alpha\pm i\omega$  for some $\alpha,\, 
\omega\in\mathbb{R},\, \omega \not=0$, then
\begin{equation}\label{E: EllipticC}\notag
e^{\pm iA}=e^{\pm i\alpha}[(\cosh{\omega }) I\pm i(\sinh{\omega})J],
\end{equation}
where 
$J=\frac{1}{\omega}(A-\alpha I)$. 

	\textbf{Case 3}: If $\lambda_1\not=\lambda_2$ (both real), then  
\begin{equation}\label{E: HyperbolicC}\notag
e^{\pm iA}=e^{\pm i\alpha }\{\,(\cos{\beta}) I\pm i(\sin{\beta})J\;\},
\end{equation}
where 
$\alpha=\frac{\lambda_1+\lambda_2}{2},\; \beta=\frac{\lambda_1-\lambda_2}{2}$
and $J=\frac{1}{\beta}(A-\alpha I)$. 
\end{corollary}
\Proof These formulas follow directly from Theorem~\ref{T: Main 1} for $t=i$. $\blacktriangle$
\begin{corollary}[Trigonometric Functions]\label{C: Sin} Let $A$ be a $2\times2$-matrix with real entries
and let
$\lambda_1$ and $\lambda_2$ be the characteristic roots of $A$. 

		\textbf{Case 1}: If $\lambda_1=\lambda_2$ (real), then 
\begin{equation}\label{E: ParabolicSin}\notag
\cos{A}=(\cos{\lambda_0})I-(\sin{\lambda_0})N,\quad \sin{A}=(\sin{\lambda_0})I+(\cos{\lambda_0})N, 
\end{equation}
where
$\lambda_0=\lambda_1=\lambda_2$ and $N=A-\lambda_0 I$.

		\textbf{Case 2}: If $\lambda_{1,2}=\alpha\pm i\omega$  for some $\alpha,\, 
\omega\in\mathbb{R},\, \omega \not=0$, then 
\begin{align}\label{E: EllipticSin}\notag
&\cos{A}=(\cos{\alpha})(\cosh{\omega})I-(\sin{\alpha})(\sinh{\omega})J,\\\notag
&\sin{A}=(\sin{\alpha})(\cosh{\omega})I+(\cos{\alpha})(\sinh{\omega})J,
\end{align}
where $J=\frac{1}{\omega}(A-\alpha I)$.

	\textbf{Case 3}: If $\lambda_1\not=\lambda_2$ (both real), then 
\begin{align}\notag
&\cos{A}=(\cos{\alpha})(\cos{\beta})I-(\sin{\alpha})(\sin{\beta})J,\\\notag
&\sin{A}=(\sin{\alpha})(\cos{\beta})I+(\cos{\alpha})(\sin{\beta})J.
\end{align}
where 
$\alpha=\frac{\lambda_1+\lambda_2}{2},\; \beta=\frac{\lambda_1-\lambda_2}{2}$
and $J=\frac{1}{\beta}(A-\alpha I)$.
\end{corollary}
\Proof As before, the result follows directly from Corollary~\ref{C: Complex Exponent} and the formulas
$\cos{x}=\frac{1}{2}(e^{ix}+e^{-ix})$ and $\sin{x}=\frac{1}{2i}(e^{ix}-e^{-ix})$.
$\blacktriangle$
\section{Exponential of a $3\times3$ Matrix}\label{S: Exponential of a 3times3-Matrix}
In this section we extend the result about $e^A$, presented in Section~\ref{S: Exponential of a
2times2-Matrix}, for $3\times 3$ matrices $A$. We should notice that: (a) the hardest (and the most time
consuming) parts of our formulas for $e^{A}$ are the calculating the squares of some matrices; we need
neither to solve linear systems of equations, nor to invert matrices: (b) our framework is always
$\mathbb{R}^3$ (never $\mathbb{C}^3$). For that reason we believe that our method has advantages over the
most conventional methods based on the choice of a particular canonical basis with transition matrix $T$ and
calculating its inverse
$T^{-1}$. To support this point we present several examples for
solving linear systems of ordinary differential equations in three unknowns. 

	In what follows we shall use the fact that if $A$ is a matrix and
$B=A-aI$ for some real $a$, then $f_B(\lambda)=f_A(\lambda+a)$. 

\begin{theorem}\label{T: Main 2} Let $A$ be a $3\times 3$ matrix with real entries and let $\lambda_1,
\lambda_2$ and
$\lambda_3$ be the characteristic roots of $A$. Then:

	\textbf{Case 1:} If $\lambda_1=\lambda_2=\lambda_3=\lambda_0$, then $N=A-\lambda_0 I$ is a nilpotent
matrix such that $N^3=0$. Consequently, for every $t\in\mathbb{R}$ we have
\begin{equation}\label{E: Case 1}
e^{tA}=e^{\lambda_0t}\left(I+tN+\frac{t^2}{2}N^2\right).
\end{equation}

	\textbf{Case 2:} If $\lambda_1=\lambda_2\not=\lambda_3$, then $A=\lambda_0 I+(\lambda_3-\lambda_0)P+N$,
where $\lambda_0=\lambda_1=\lambda_2$, 
\[
P=\frac{1}{(\lambda_3-\lambda_0)^2}\left(A-\lambda_0I\right)^2 \text{\;and\;\,} N=A-\lambda_0
I-(\lambda_3-\lambda_0)P.
\]
Also, we have $PN=NP=0,\quad P^2=P\quad
\text{and}\quad N^2=0$. Consequently, for every $t\in\mathbb{R}$ we have
\begin{equation}\label{E: Case 2}
e^{tA}=e^{\lambda_0t}(I-P+tN)+e^{\lambda_3t}P.
\end{equation}

		\textbf{Case 3:} Let $\lambda_{1,2}=\alpha\pm i\omega$ for some $\alpha,
\omega\in\mathbb{R},\, \omega\not=0$. Then $A=\alpha I+(\lambda_3-\alpha)P+\omega J$,
where
\begin{equation}\notag
P=\frac{(A-\alpha I)^2+\omega^2 I}{(\lambda_3-\alpha)^2+\omega^2},\quad\quad
J=\frac{1}{\omega}\left[A-\alpha I-(\lambda_3-\alpha)P\right].
\end{equation}
Also, we have $PJ=JP=0,\; P^2=P,\; J^{2n}=(-1)^n(I-P),\; n=1, 2,\dots,\; J^{2n+1}=(-1)^n J,\;
n=0,1,2,\dots$. Consequently, for every $t\in\mathbb{R}$ we have
\begin{equation}\label{E: Case 3}
e^{tA}=e^{\alpha t}\left[\,(\cos{\omega
t})(I-P)+(\sin{\omega t})J\,\right]+e^{\lambda_3 t}P,
\end{equation}

	\textbf{Case 4:} If $\lambda_1, \lambda_2, \lambda_3$ are three distinct reals, then $A=\alpha
I+(\lambda_3-\alpha)P+\beta J$, where $\alpha =(\lambda_1+\lambda_2)/2,\;  \beta =(\lambda_1-\lambda_2)/2$,
and
\begin{equation}\notag
P=\frac{(A-\alpha I)^2-\beta^2 I}{(\lambda_3-\alpha)^2-\beta^2},\quad\quad
J=\frac{1}{\beta}\left[A-\alpha I-(\lambda_3-\alpha)P\right].
\end{equation}
Also, $PJ=JP=0,\; P^2=P,\; J^{2n}=I-P,\; n=1, 2\dots,\; J^{2n+1}=J,\; n=0, 1, 2\dots$.
Consequently, for every $t\in\mathbb{R}$ we have

\begin{align}\label{E: Case 4}
e^{tA}=&e^{\alpha t}\left[\,(\cosh{\beta
t})(I-P)+(\sinh{\beta t})J\,\right]+e^{\lambda_3 t}P\overset{or}=\\\notag
&\overset{or}=\frac{1}{2}e^{\lambda_1 t}(I+J-P)+
\frac{1}{2}e^{\lambda_2 t}(I-J-P)+e^{\lambda_3 t}P.
\end{align}
\end{theorem}

\Proof Case 1: We have $f_A(\lambda)=(\lambda-\lambda_0)^3$ and
$N^3=(A-\lambda_0I)^3=f_A(A)=0$, by the Cayley-Hamilton
theorem. Next,
$e^{tA}=e^{t(\lambda_0I+N)}=e^{t\lambda_0I}\,
e^{tN}=e^{\lambda_0t}I(I+tN+(t^2/2)N^2+0+\dots)=e^{\lambda_0t}(I+tN+(t^2/2)N^2)$, as
required. 

	Case 2: Denote $\lambda_3-\lambda_0=b,\; A-\lambda_0I=B$ and observe that $P=\frac{1}{b^2}B^2$ and
$N=B-bP=-\frac{1}{b}B(B-bI)$. Notice that $f_A(\lambda)=(\lambda-\lambda_0)^2(\lambda-\lambda_3)$ implying,
by translation,
$f_B(\lambda)=\lambda^2(\lambda-b)$. It follows $B^2(B-bI)=f_B(B)=0$ (by the Cayley-Hamilton
theorem), i.e. $B^3=bB^2$ and $B^4=b^2B^2$. After these preliminary evaluations we calculate
$P^2=B^4/b^4=b^2B^2/b^4=P$, as required. Next, we have $NP=PN=-\frac{1}{b^3}B^2B(B-bI)=
-\frac{1}{b^3}Bf_B(B)=0$, and
$N^2=-\frac{1}{b^2}B^2(B-bI)^2=\frac{1}{b^2}f_B(B)(B-bI)=0$, as required. Finally,
\begin{align}\notag
&e^{tA}=e^{t[\lambda_0I+bP+N]}=e^{\lambda_0tI}\,e^{btP}\,e^{tN}=\\\notag
&=e^{\lambda_0t}I\left(I+P\sum_{n=1}^\infty\frac{t^nb^n}{n!}\right)(I+tN)=
e^{\lambda_0t}[I-P+e^{bt}P](I+tN)=\\\notag
&=e^{\lambda_0t}(I-P+tN)+e^{\lambda_3t}P,
\end{align}
as required. 

	Case 3: We denote $A-\alpha I=B$ and $\lambda_3-\alpha=b$ and observe that
$P=\frac{1}{b^2+\omega^2}(B^2+\omega^2I)$ and $J=\frac{1}{\omega}(B-bP)$. We have
$f_A(\lambda)=[(\lambda-\alpha)^2+\omega^2](\lambda-\lambda_3)$ which implies
$f_B(\lambda)=(\lambda^2+\omega^2)(\lambda-b)$. By the Cayley-Hamilton
theorem, it follows $(B^2+\omega^2I)(B-bI)=0$ implying also $B^3+\omega^2B=b(B^2+\omega^2I)$
and $B^4+\omega^2B^2=b^2(B^2+\omega^2I)$. Next we calculate
$PB=BP=\frac{1}{b^2+\omega^2}(B^3+\omega^2B)=
\frac{b}{b^2+\omega^2}(B^2+\omega^2I)=bP$. It follows $B^2P=b^2P$ which
helps us to show that $P$ is a projection, i.e.
$P^2=\frac{1}{b^2+\omega^2}(B^2+\omega^2I)P=\frac{1}{b^2+\omega^2}(B^2P+\omega^2P)=
\frac{1}{b^2+\omega^2}(b^2P+\omega^2P)=P$ and also
$PJ=JP=\frac{1}{\omega}(B-bP)P=\frac{1}{\omega}(bP-bP)=0$. Finally, we have
$J^2=\frac{1}{\omega^2}(B-bP)^2=\frac{1}{\omega^2}(B^2-2bBP+b^2P)=\frac{1}{\omega^2}(B^2-b^2P)
=\frac{1}{\omega^2}[(b^2+\omega^2)P-\omega^2I-b^2P]=-(I-P)$, as required. The rest of the formulas for the
powers of $J$ follow immediately. Next, we calculate
$e^{\omega tJ}$: 
\begin{align}\notag
&e^{\omega tJ}=\sum_{n=0}^\infty\frac{(\omega t)^n}{n!}J^n=I+\sum_{n=1}^\infty\frac{(\omega
t)^{2n}}{(2n)!}J^{2n}+\sum_{n=0}^\infty\frac{(\omega
t)^{2n+1}}{(2n+1)!}J^{2n+1}=\\\notag
&=P+ (I-P)+\sum_{n=1}^\infty\frac{(-1)^n(\omega
t)^{2n}}{(2n)!}(I-P)+\sum_{n=0}^\infty\frac{(-1)^n(\omega t)^{2n+1}}{(2n+1)!}J=\\\notag
&=P+(\cos{\omega t})(I-P)+(\sin{\omega t})J.
\end{align}
Finally, we calculate 
\begin{align} \notag
&e^{tA}=e^{t[\alpha I+bP+\omega J]}=
e^{\alpha tI}\,e^{btP}\,e^{\omega tJ}=
e^{\alpha t}[(I-P)+e^{bt}P]\times\\\notag
&\times [P+(\cos{\omega t})(I-P)+(\sin{\omega t})J].
\end{align}
The last leads to formula (\ref{E: Case 3}) after standard manipulations. 

	Case 4: We denote, as before, $A-\alpha I=B$ and $\lambda_3-\alpha=b$ and observe that
$P=\frac{1}{b^2-\beta^2}(B^2-\beta^2I)$ and $J=\frac{1}{\beta}(B-bP)$. Next, we have
$f_A(\lambda)=(\lambda-\lambda_1)(\lambda-\lambda_2)(\lambda-\lambda_3)$ which implies, by translation,
$f_B(\lambda)=(\lambda^2-\beta^2)(\lambda-b)$. It follows
$(B^2-\beta^2I)(B-bI)=0$, by the Cayley-Hamilton theorem. Thus
$B^3-\beta^2B=b(B^2-\beta^2I)$ and $B^4-\beta^2B^2=b^2(B^2-\beta^2I)$. Now the relations between $P$ and $J$
follow exactly as in Case 3. Next we calculate:
\begin{align} \notag
&e^{tA}=e^{t[\alpha I+bP+\beta J]}=e^{\alpha tI}\,e^{btP}\,e^{\beta t
J}=e^{\alpha t}[I-P+e^{bt}P]\times\\\notag
&\times\left[I+\sum_{n=1}^\infty\frac{(\beta t)^{2n}}{(2n)!}(I-P)+\sum_{n=0}^\infty\frac{(\beta
t)^{2n+1}}{(2n+1)!}J\right]=[e^{\alpha t}(I-P)+e^{\lambda_3t}P\,]\times\\\notag
&\times\left[P+(\cosh{\beta t})(I-P)+(\sinh{\beta t})J\right]=e^{\alpha t}\left[(\cosh{\beta
t})(I-P)+(\sinh{\beta t})J\right]+\\\notag
&e^{\lambda_3 t}P,
\end{align} 
as required. $\blacktriangle$

	The formulas (\ref{E: Case 1})-(\ref{E: Case 4}) show that we do not need the eigenvectors of $A$ in order
to calculate $e^{tA}$. However, the eigenvectors of $A$ as well as canonical forms and canonical bases (if
needed for classification, graphing, etc.) can be easily extracted from the matrices
$N, P$ and $J$. The next corollary follows easily from the above theorem and we leave the proof to the
reader. 
\begin{corollary}[Eigenvectors, Canonical Forms and Bases]\label{C: Eigenvectors, Canonical Forms and
Bases 3by3} Under the notation of the previous theorem we have the following:

	\textbf{Case 1:} If $N^2\not=O$, then $E_{\lambda=\lambda_0}(A)= {\rm Im}(N^2)$. In particular,
either of the non-zero columns of $N^2$ is an eigenvector of $A$. We
have
$A=T\begin{pmatrix}\lambda_0&1&0\\0&\lambda_0&1\\0&0&\lambda_0\end{pmatrix}T^{-1},$ where
$T=[\, N^2\vec{x}| N\vec{x}\,|\,\vec{x}\,]$ and
$\vec{x}\in\mathbb{R}^3$, $N^2\vec{x}\not=\vec{0}$. If $N=O$, then $A=\lambda_0 I$.  If $N^2=O$ and
$N\not=O$, then $E_{\lambda=\lambda_0}(A)\supsetneqq{\rm Im}(N)$. To complete a Jordan basis we can select
an additional eigenvector by easy inspection of the columns of $N$. 

	\textbf{Case 2:} If $N\not=O$, then $E_{\lambda=\lambda_0}(A)= {\rm Im}(N)$ and $E_{\lambda=\lambda_3}(A)=
{\rm Im}(P)$. In particular, either of the non-zero columns of $N$ and $P$ is an eigenvector of $A$. We
have
$A=T\begin{pmatrix}\lambda_0&1&0\\0&\lambda_0&1\\0&0&\lambda_3\end{pmatrix}T^{-1},$ where

\[
T=[\, (A-\lambda_0I)(A-\lambda_3I)\vec{x}\,|(A-\lambda_3I)\vec{x}|\,P\vec{y}\,],
\]
and $\vec{x}\in\mathbb{R}^3$, $(A-\lambda_0I)(A-\lambda_3I)\vec{x}\not=\vec{0}$ and $P\vec{y}\not=\vec{0}$.
If
$N=O$, then
$A={\rm diag}(\lambda_0, \lambda_0, \lambda_3)$.  

	\textbf{Case 3:} The only eigenspace in $\mathbb{R}^3$ is $E_{\lambda=\lambda_3}={\rm Im}(P)$. Here we
have $A=T\begin{pmatrix}\alpha & \omega&0\\-\omega &\alpha&0\\0&0&\lambda_3\end{pmatrix}
T^{-1}$, where $T=[\,J^2\vec{x}\,|J\vec{x}\,|\,P\vec{y}\, ]$ and $\vec{x}, \vec{y}\in\mathbb{R}^3,
\vec{x}\not=\vec{0}$, $P\vec{y}\not=\vec{0}$. Alternatively, in 
$\mathbb{C}^3$ the matrix $A$ is diagonalizable with eigenvalues $\alpha\pm i\omega,\, \lambda_3$ and
eigenvectors $(J\mp iJ^2)\vec{x},\, \vec{x}\not=\vec{0}$, $P\vec{y}\not=\vec{0}$, respectively.

	\textbf{Case 4:} We have
$E_{\lambda=\lambda_{1,2}}(A)= {\rm Im}(J\pm J^2)$, respectively, and $E_{\lambda=\lambda_3}(A)= {\rm
Im}(P)$. Thus the matrix $A$ is diagonalizable in the basis consisting of the eigenvectors listed above.
\end{corollary}

	Here are several examples.
\begin{example} Let $\vec{x}^{\,\prime}(t)=A\vec{x}(t)$, where 
$A=\begin{pmatrix}2 & -1 & 2\\5 & -3 &3\\-1&0&-2\end{pmatrix}$, be a system of three first order
homogeneous ordinary differential equations with constant coefficients with initial conditions
$\vec{x}(\vec{0})=(1, 0, 2)$. We apply the formula for the solution
$\vec{x}(t)=e^{tA}\,\vec{x}_0$. It remains to calculate $e^{tA}$. We leave to the reader to check that
$\lambda_1=\lambda_2=\lambda_3=-1$ (Case 1 in Theorem~\ref{T: Main 2}). We calculate 
\[
N=A+I=\begin{pmatrix}3 & -1 & 2\\5 & -2 &3\\-1&0&-1\end{pmatrix},\; N^2=\begin{pmatrix}2 & -1 & 1\\2 &
-1 &1\\-2&1&-1\end{pmatrix},\; N^3=O.
\]
The formula (\ref{E: Case 1}) gives
\begin{equation}\notag
e^{tA}=e^{\lambda_0t}\left(I+tN+\frac{t^2}{2}N^2\right)=e^{-t}\begin{pmatrix}1+3t+t^2, & -t-t^2/2, &
2t+t^2/2\\5t+t^2, & 1-2t-t^2/2, &3t+t^2/2\\-t-t^2,&t^2/2,&1-t-t^2/2\end{pmatrix}.
\end{equation}
The solution of the system is
$\vec{x}(t)=e^{tA}\begin{pmatrix}1\\0\\2\end{pmatrix}=
e^{-t}\begin{pmatrix}1+7t+2t^2\\11t+2t^2\\2-3t-2t^2
\end{pmatrix}$.
\end{example}

\begin{remark}[Jordan Basis and Eigenvectors] We do not need a Jordan basis (and the eigenvectors) of
$A$ to calculate
$e^{tA}$. However, a Jordan basis for $A$, if needed, can be easily extracted from
the above calculations with the help of Corollary~\ref{C: Eigenvectors, Canonical Forms and
Bases 3by3}. For example, $\{N^2\vec{e_3},\, N\vec{e_3},\,
\vec{e_3}\}$ (where $\vec{e_3}$ is the third vector in the standard basis in $\mathbb{R}^3$), forms a Jordan
basis for
$A$. Also, we have $E_{\lambda=-1}={\rm Im}(N^2)$. In particular, either of the columns of $N^2$ is an
eigenvector of $A$ for $\lambda=-1$. We leave to the reader to check that
$A=TCT^{-1}$, where $C=\begin{pmatrix}-1 & 1 & 0\\0 & -1 &1\\0&0&-1\end{pmatrix}$ and
$T=[N^2\vec{e_3}\,|\, N\vec{e_3}\,|\, \vec{e_3}]=\begin{pmatrix}1 & 2 & 0\\1 & 3
&0\\-1&-1&1\end{pmatrix}$. The columns of $T$ forms a Jordan basis for $A$.
\end{remark}

\begin{example} Let $\vec{x}^{\,\prime}(t)=A\vec{x}(t)$, where 
$A=\begin{pmatrix}1 & -3 & 4\\4 & -7 &8\\6&-7&7\end{pmatrix}$, be a system of three first order
homogeneous ordinary differential equations with constant coefficients with initial conditions
$\vec{x}(\vec{0})=(1, 0, 1)$. We leave to the reader to check that
$\lambda_1=\lambda_2=-1=\lambda_0$ and $\lambda_3=3$ (Case 2 in Theorem~\ref{T: Main 2}). We calculate 
\begin{align}\notag
&A+I=\begin{pmatrix}2 & -3 & 4\\4 & -6 &8\\6&-7&8\end{pmatrix},\quad
P=\left(\frac{A+I}{4}\right)^2=\begin{pmatrix}1 & -1 & 1\\2 & -2 &2\\2&-2&2\end{pmatrix},\\\notag
& N=A+I-4P=\begin{pmatrix}2 & -3 & 4\\4 & -6 &8\\6&-7&8\end{pmatrix}-4\begin{pmatrix}1 & -1 & 1\\2 & -2
&2\\2&-2&2\end{pmatrix}=\begin{pmatrix}-2 & 1 & 0\\-4 & 2 &0\\-2&1&0\end{pmatrix}.
\end{align}
Notice that $N^2=O$. The formula (\ref{E: Case 2}) gives
\begin{align}\notag
e^{tA}=&e^{-t}\begin{pmatrix}1-2t & t & 0\\-4t & 1+2t
&0\\-2t&t&1\end{pmatrix}+(e^{3t}-e^{-t})\begin{pmatrix}1 & -1 & 1\\2 & -2
&2\\2&-2&2\end{pmatrix}=\\\notag
&=\begin{pmatrix}-2te^{-t}+e^{3t}, & (1+t)e^{-t}-e^{3t}, & -e^{-t}+e^{3t}\\(-2-4t)e^{-t}+2e^{3t}, &
(3+2t)e^{-t}-2e^{3t},
&-2e^{-t}+2e^{3t}\\-4te^{-t}+2e^{3t},&(2+t)e^{-t}-2e^{3t},&-2e^{-t}+2e^{3t}\end{pmatrix}.
\end{align}

For the solution we have
$\vec{x}(t)=e^{tA}\begin{pmatrix}1\\0\\1\end{pmatrix}=
\begin{pmatrix}(-1-2t)e^{-t}+2e^{3t}\\(-4-4t)e^{-t}+4e^{3t}\\(-2-4t)e^{-t}+4e^{3t}
\end{pmatrix}$.
\end{example}
\begin{remark}[Jordan Basis and Eigenvectors] We just calculated $e^{tA}$ without the Jordan basis and
eigenvectors of $A$. They, however (if needed),
can be easily extracted from the above calculations with the help of 
Corollary~\ref{C: Eigenvectors, Canonical Forms and
Bases 3by3}. For example,
$E_{\lambda=-1}={\rm Im}(N)$ and $E_{\lambda=3}={\rm Im}(P)$. In particular, the first and second
columns of $N$ are eigenvectors of $A$ for 
$\lambda=-1$ and either of the columns of $P$ is an eigenvector of $A$  and $\lambda=3$. We leave to the
reader to check that
$A=TCT^{-1}$, where $C=\begin{pmatrix}-1 & 1 & 0\\0 & -1 &0\\0&0&3\end{pmatrix}$ and
$T=[\, (A-\lambda_0I)(A-\lambda_3I)\vec{x}\,|(A-\lambda_3I)\vec{x}|\,P\vec{y}\,]=\begin{pmatrix}8 & -1 &
1\\16 & 2 &2\\8&3&2\end{pmatrix}$ for $\vec{x}=\vec{y}=\vec{e_1}$. The columns of
$T$ forms a Jordan basis for $A$.
\end{remark}
\begin{example} Let $\vec{x}^{\,\prime}(t)=A\vec{x}(t)$, where 
$A=\begin{pmatrix}1 &1 & -1\\0 &3 &0\\1&0&1\end{pmatrix}$, be a system of three first order
homogeneous ordinary differential equations with constant coefficients with initial conditions
$\vec{x}(\vec{0})=\vec{e_2}$, where $\vec{e_2}$ is the second vector in the standard basis in
$\mathbb{R}^3$. We have $\lambda_{1,2}=1\pm i$, i.e. $\gamma=\omega=1$, and $\lambda_3=3$ (Case 3 in
Theorem~\ref{T: Main 2}). We calculate 
\[
P=\frac{(A-\gamma I)^2+\omega^2I}{(\lambda_3-\gamma)^2+\omega^2}=
\frac{1}{5}\begin{pmatrix}0 & 2 &0\\0 &5 &0\\0&1&0\end{pmatrix}, 
\]
\begin{align}\notag
&J=\frac{1}{\omega}\left[A-\gamma I-(\lambda_3-\gamma)P\right]=
\frac{1}{5}\begin{pmatrix}0 & 1 &-5\\0 &0&0\\5&-2&0\end{pmatrix},\\\notag 
& e^{tA}=(e^{\lambda_3 t}-e^{\gamma t}\cos{\omega t})P+e^{\gamma t}[(\cos{\omega t}) I+(\sin{\omega
t})J]=\\\notag
&=\frac{1}{5}\begin{pmatrix}5e^t\cos{t}, & 2e^{3t}-e^t(2\cos{t}-\sin{t}),&-5e^t\sin{t}\\
0,&5e^{3t},&0\\5e^t\sin{t},&e^{3t}-e^t(\cos{t}+2\sin{t}),&5e^t\cos{t}\end{pmatrix}.
\end{align}
	For the solution we have
$\vec{x}(t)=e^{tA}\begin{pmatrix}0\\1\\0\end{pmatrix}=
\frac{1}{5}\begin{pmatrix}2e^{3t}-e^t(2\cos{t}-\sin{t})\\5e^{3t}\\e^{3t}-e^t(\cos{t}+2\sin{t})
\end{pmatrix}$.
\end{example}
\begin{remark}[Canonical Forms and Eigenvectors] As before, we just calculated $e^{tA}$ without the Jordan
basis and eigenvectors of $A$. They, however (if needed),
can be easily extracted from the above calculations with the help of 
Corollary~\ref{C: Eigenvectors, Canonical Forms and
Bases 3by3}. In this case $E_{\lambda=\lambda_3}={\rm Im}(P)$ and
we have $A=TCT^{-1}$, where $C=\begin{pmatrix}1 & 1 & 0\\-1 & 1 &0\\0&0&3\end{pmatrix}$ and
$T=[\, J^2\vec{x}\,|J\vec{x}|\,P\vec{y}\,]=\begin{pmatrix}-1 & 0 & 2\\0 & 0 &5\\0&1&1\end{pmatrix}$ for
$\vec{x}=\vec{e_1}$ and $\vec{y}=5\vec{e_2}$. The columns of
$T$ form a canonical basis for $A$.
\end{remark}



\end{document}